\documentclass[11pt,reqno]{amsart}

\usepackage{amsmath,amssymb}
\usepackage{fancyhdr,version}
\usepackage{fullpage}
\usepackage{longtable}
\usepackage{color}

\theoremstyle{plain} \newtheorem{theorem}{Theorem}

\theoremstyle{plain} 

\theoremstyle{plain} \newtheorem{lemma}{Lemma}

\theoremstyle{definition} 

\theoremstyle{plain}

\newcommand{\D}{\mathbb{D}}

\newcommand{\Sch}{\mathcal{S}}

\newcommand{\sg}{\sigma}

\newcommand{\Hess}{\operatorname{Hess}}

\renewcommand{\Re}{\operatorname{Re}}

\newcommand{\sff}{\ensuremath{I\!I}}

\newcommand{\g}{\mathbf{g}}

\title{Ahlfors-Weill Extensions for a Class of Minimal Surfaces}

\thanks{The authors were supported in part by FONDECYT Grant \# 1071019.}

\author{M. Chuaqui \and P. Duren \and B. Osgood}

\address{P. Universidad Cat\'olica
de Chile}
\email{mchuaqui@mat.puc.cl}
\address{University of Michigan}
\email{duren@umich.edu}
\address{Stanford University}
\email{osgood@stanford.edu}



\begin{document}

\maketitle

\begin{abstract}
The Ahlfors-Weill extension of a conformal mapping of the disk is generalized to the lift of a harmonic mapping of the disk to a minimal surface, producing homeomorphic and quasiconformal extensions. The extension is obtained by a reflection across the boundary of the surface using a family of Euclidean circles orthogonal to the surface. This gives a geometric generalization of the Ahlfors-Weill formula and extends the minimal surface.  Thus one obtains a homeomorphism of $\overline{\mathbb{C}}$ onto a toplological sphere in $\overline{\mathbb{R}^3} = \mathbb{R}^3 \cup \{\infty\}$ that is  real-analytic off the boundary. The hypotheses involve bounds on a generalized Schwarzian derivative for harmonic mappings in term of the hyperbolic metric of the disk and the Gaussian curvature of the minimal surface. Hyperbolic convexity plays a crucial role.

\end{abstract}

\section{Introduction \label{section:intro}}

 If $f$ is an analytic, locally injective function its Schwarzian derivative is
\[
Sf = \left(\frac{f''}{f'}\right)'-\frac{1}{2}\left(\frac{f''}{f'}\right)^2.
\]
We owe to Nehari \cite{nehari:schlicht} the discovery that the size of the Schwarzian derivative of an analytic function is related to its injectivity, and to Ahlfors and Weill \cite{ahlfors-weill:extension} the discovery of an  allied, stronger phenomenon of quasiconformal extension of the function. We  state the combined results as follows:
\begin{theorem} \label{theorem:plane-case}
Let $f$ be analytic and locally injective in the unit disk, $\D$.
\begin{list}{}{\setlength\leftmargin{.25in}}
\item[(a)] If
\begin{equation} \label{eq:nehari-2}
|Sf(z)| \le \frac{2}{(1-|z|^2)^2}, \quad z \in \D,
\end{equation}
then $f$ is injective in $\D$.
\item[(b)] If for some $t<1$
\begin{equation} \label{eq:aw-2t}
|Sf(z)| \le\frac{2t}{(1-|z|^2)^2}, \quad z \in \D,
\end{equation}
then $f$ has a $\frac{1+t}{1-t}$-quasiconformal extension to $\overline{\mathbb{C}}$.
\end{list}
\end{theorem}
A remarkable aspect of Ahlfors and Weill's theorem is the explicit formula they give for the quasiconformal extension. They need the stronger inequality \eqref{eq:aw-2t} to show, first of all, that  the extended mapping has a positive Jacobian and is hence a local homeomorphism. Global injectivity then follows from the monodromy theorem and quasiconformality from a calculation of the dilatation. The topological argument cannot get started without \eqref{eq:aw-2t}, but a different approach in \cite{co:aw} shows that the same formula still provides a homeomorphic extension even when $f$ satisfies the weaker inequality \eqref{eq:nehari-2} and   $f(\D)$ is a Jordan domain. As to the latter requirement, if $f$ satisfies \eqref{eq:nehari-2} then $f(\D)$ fails to be a Jordan domain only when $f(\D)$ is a parallel strip or the image of a parallel strip under a M\"obius transformation, as shown by Gehring and Pommerenke \cite{gp:nehari}.

 \smallskip

In earlier work we introduced a Schwarzian derivative for plane harmonic mappings and we established an injectivity criterion analogous to \eqref{eq:nehari-2} for the Weierstrass-Enneper lift of a harmonic mapping of $\D$  to a minimal surface. In this paper we show that injective and quasiconformal extensions also obtain in this more general setting under  conditions analogous to \eqref{eq:nehari-2} and \eqref{eq:aw-2t}, respectively.  The construction is a geometric generalization of the Ahlfors-Weill formula and extends the minimal surface.  Thus one obtains a homeomorphism of $\overline{\mathbb{C}}$ onto a toplological sphere in $\overline{\mathbb{R}^3} = \mathbb{R}^3 \cup \{\infty\}$ that is actually real-analytic off $\partial\D$. Precise statements require some additional preparation, and for more background and details we refer to \cite{cdo:injective-lift}. 

Let $f \colon \D \longrightarrow \mathbb{C}$ be a harmonic mapping. As is customary we write $f=h+\bar{g}$, where $g$ and $h$ are analytic. We assume that $f$ is locally injective and that the dilatation $\omega = g'/h'$ is the square of a meromorphic function on $\D$. Under these assumptions there is a lift $\widetilde{f}\colon \D \longrightarrow \Sigma$, the Weierstrass-Enneper lift, onto a minimal surface $\Sigma \subset \mathbb{R}^3$. Furthermore, $\widetilde{f}$ is a conformal mapping of $\D$ to $\Sigma$, each with its Euclidean metric. We let $\mathbf{g}_0$ denote the Euclidean metric on $\mathbb{R}^3$, or the induced Euclidean metric on $\Sigma$.  The pullback of $\mathbf{g}_0$ is a conformal metric on $\D$:
\[
e^{2\sg}|dz|^2 = \widetilde{f}^*(\mathbf{g}_0) \quad \text{where} \quad e^{\sg} = |h'|+|g'|.
\]
In terms of $\sigma$, the Gauss curvature of $\Sigma$ at a point $\widetilde{f}(z)$ is
\[
K(\widetilde{f}(z)) = -e^{-2\sg(z)}\Delta \sg(z). 
\]
For a minimal surface the curvature is $\le 0$. The Schwarzian of $f$ (or of $\widetilde{f}$)  is
\begin{equation} \label{eq:harmonic-schwarzian}
\Sch f =  2(\sg_{zz}-\sg_z^2).
\end{equation}
This becomes the usual Schwarzian when $f$ is a analytic, in which case $\sg=\log|f'|$.

Much of our work will go into defining an injective, continuous reflection of $\Sigma$ across its boundary, $R \colon {\Sigma} \longrightarrow \Sigma^* \subset \overline{\mathbb{R}^3}$,  with which we will extend $\widetilde{f}$ to
\[
\widetilde{F}(z)=
\begin{cases}
\widetilde{f}(z), & \quad z \in\overline{\D},\\
R(\widetilde{f}(1/\bar{z})), & \quad z \in \overline{\mathbb{C}}\setminus \overline{\D}.
\end{cases}
\]
The analysis will include a discussion of boundary values. 

We state our results in parallel to Theorem \ref{theorem:plane-case}, including the homeomorphic extension for the first part: 
\begin{theorem} \label{theorem:extension}
Let $f$ be harmonic and locally injective in $\D$ with lift $\widetilde{f} \colon \D \longrightarrow \Sigma$.
\begin{list}{}{\setlength\leftmargin{.25in}}
\item[(a)] If
\begin{equation} \label{eq:nehari-harmonic}
|\mathcal{S}f(z)| + e^{2\sigma(z)} |K(\widetilde{f}(z))| \leq \frac{2}{(1-|z|^2)^2}\,,
\quad z\in\Bbb D\,, 
\end{equation}
then $\widetilde{f}$ is injective in $\D$. If $\widetilde{f}(\partial\mathbb{D})$ is a Jordan curve then $\widetilde{F}$ is a continuous, injective extension to $\overline{\mathbb{C}}$.
\item[(b)] If for some $t<1$
\begin{equation} \label{eq:aw-harmonic}
|\mathcal{S}f(z)| + e^{2\sigma(z)} |K(\widetilde{f}(z))| \leq \frac{2t}{(1-|z|^2)^2}\,,
\quad z\in\Bbb D\,, 
\end{equation}
and if for some constant $C$
\begin{equation} \label{eq:sigma-condition}
\|\nabla\sg(z)\| \le \frac{C}{1-|z|^2},\quad z \in \D,
\end{equation}
then  $\widetilde{F}$ is a quasiconformal extension to $\overline{\mathbb{C}}$ with a bound depending only on $t$ and $C$.
\end{list}
\end{theorem}

The injectivity in part (a)  was proved in \cite{cdo:injective-lift}  in even greater generality, so the point here is the extension.   It was also proved in \cite{cdo:injective-lift} that if $f$ satisfies \eqref{eq:nehari-harmonic} then $f$ and $\widetilde{f}$ have spherically continuous extensions to $\partial\D$. Furthermore, we know exactly when $\widetilde{f}(\partial\D)$ fails to be a simple closed curve in $\mathbb{R}^3$, namely when  $\widetilde{f}$ maps $\overline{\D}$ into a catenoid  and $\partial\Sigma$ is pinched by a \emph{Euclidean circle} on the surface. More precisely, there is a Euclidean circle $C$ on $\overline{\Sigma}$ and a point $P\in C$ with $\widetilde{f}(\zeta_1) = P=\widetilde{f}(\zeta_2)$ for a pair of points $\zeta_1,\zeta_2 \in \partial\D$.  Equality holds in \eqref{eq:nehari-harmonic} along $\widetilde{f}^{-1}(C \setminus\{P\})$, and because of this a  function satisfying the stronger inequality \eqref{eq:aw-harmonic} is always injective on $\partial\D$. 

Independent of its connection with injectivity, an enduring source of interest in the analytic Schwarzian stems from its invariance properties under M\"obius transformations:  if $T(z) = (az+b)/(cz+d)$ then
\begin{equation} \label{eq:Tf}
S(T\circ f) = Sf \quad \text{and} \quad S(f\circ T) = ((Sf) \circ T)(T')^2.
\end{equation}
For harmonic mappings and the harmonic Schwarzian the former equation does not apply since $T\circ f$ is generally not harmonic.  However, the latter equation continues to hold. As a consequence of this and Schwarz's Lemma,  if a harmonic mapping $f$ satisfies \eqref{eq:nehari-harmonic} or \eqref{eq:aw-harmonic} and if $T$ is a M\"obius transformation of $\D$ onto itself, then $f\circ T$ also satisfies the inequalities. The equations \eqref{eq:Tf} are contained in the more general chain rule for the Schwarzian,
\begin{equation} \label{eq:chain-rule}
S(g\circ f) = ((Sg)\circ f)(f')^2 +Sf.
\end{equation}

By `quasiconformal' we mean that $\widetilde{F}$ satisfies
\begin{equation} \label{eq:qc-bound}
\frac{\max_{\|X\|=1}\|D_X\widetilde{F}\|}{\min_{\|X\|=1}\|D_X\widetilde{F}\|} \le A
\end{equation}
at all points in $\mathbb{C}\setminus \partial\D$ for an $A$ that depends only on $t$ and $C$. The ratio is $1$ at points in $\D$ because there $\widetilde{F}(z) = \widetilde{f}(z)$ is conformal.

The statements in Theorem \ref{theorem:extension} all reduce to their classical counterparts when $f$ is analytic, including the formula for the extension. The condition \eqref{eq:sigma-condition} becomes
\[
\left|\frac{f''(z)}{f'(z)}\right| \le \frac{C}{1-|z|^2}
\]
and this is true with $C=6$ when $f$ is injective in $\D$, in particular when $f$ satisfies \eqref{eq:nehari-2}. For harmonic mappings we must assume  \eqref{eq:sigma-condition}, but  it is a mild restriction that holds for many cases of interest, for example when $\widetilde{f}$ is bounded.

\bigskip

The proof of Theorem \ref{theorem:extension} is in several parts. In Section \ref{section:circles-reflections} we will construct the reflection and show that it is injective when $f$ satisfies \eqref{eq:nehari-harmonic}. This is supported by lemmas on convexity and critical points proved in Section \ref{section:convexity}. In Section \ref{section:extension} we show that the extension  matches up continuously along $\partial\D$, completing the proof of the first part of the theorem.  In Section \ref{section:qc-calculation} we show that the reflection, and hence the extension, is quasiconformal when $f$ satisfies the stronger inequality \eqref{eq:aw-harmonic}.

The reflection $w\mapsto  w^*$ sews a surface $\Sigma^*=R(\Sigma)$ to the minimal surface $\Sigma$ along the boundary. It would be interesting to study the geometry of $\Sigma^*$, both when $R$ is simply injective and especially when it is quasiconformal. The latter provides a class of surfaces that are quasiconformaly equivalent to a sphere, about which there is limited knowledge. We hope to return to this topic on another occasion.

\section{Three Lemmas on Convexity and Critical Points} \label{section:convexity}

In this section we borrow some results and techniques from \cite{cdo:injective-lift}, all having to do with convexity, to set the stage for constructing the reflection.

A real-valued function $u$ on $\D$ is \emph{hyperbolically convex} if
\begin{equation} \label{eq:u''>0}
(u\circ \gamma)''(s) \ge 0
\end{equation}
for all hyperbolic geodesics $\gamma(s)$ in $\D$, where $s$ is the hyperbolic arclength parameter. A special case of Theorem 4 in \cite{cdo:injective-lift} tells us that when $f$ satisfies the injectivity condition \eqref{eq:nehari-harmonic} the positive function
\begin{equation} \label{eq:u}
u_{\widetilde{f}}(z)= \frac{1}{\sqrt{(1-|z|^2)e^{\sg(z)}}}, \quad z \in \D,
\end{equation}
is hyperbolically convex. The principle is that an upper bound for the Schwarzian leads to a lower bound for the Hessian of $u_{\widetilde{f}}$, and from there to \eqref{eq:u''>0} when $u_{\widetilde{f}}$ is restricted to a geodesic. We will have some additional comments at the end of this section.

A second principle is to employ a version of the Schwarzian introduced by Ahlfors in \cite{ahlfors:schwarzian-R^n} when studying conditions such as \eqref{eq:nehari-harmonic} along curves.  Let $\varphi : (a,b)\rightarrow {\Bbb R}^n$ be of class $C^3$ with
$\varphi'(x)\neq0$. Ahlfors defined
\begin{equation}
\label{eq:S1}
S_1\varphi = \frac{\langle \varphi''',\varphi'\rangle}{|\varphi'|^2}
- 3\frac{\langle \varphi'',\varphi'\rangle^2}{|\varphi'|^4}
+ \frac32\frac{\|\varphi''\|^2}{\|\varphi'\|^2}\,, 
\end{equation}
where $\langle\cdot,\cdot\rangle$ denotes the Euclidean inner product.  If $T$ is a M\"obius transformation of $\overline{\mathbb{R}^n}$ then $S_1(T\circ \varphi) = S_1\varphi$, so this important invariance property is available.

 Whereas Ahlfors' interest was in the relation of $S_1\varphi$ to the change in  cross ratio under $\varphi$, another geometric property of $S_1\varphi$ was discovered by Chuaqui and Gevirtz in  \cite{chuaqui-gevirtz:S1}. Namely, if
 \[
 v=\|\varphi'\|
 \]
 then
 \begin{equation} \label{eq:S1-curvature}
 S_1\varphi = \left(\frac{v'}{v}\right)'-\frac{1}{2}\left(\frac{v'}{v}\right)^2 + \frac{1}{2}v^2\kappa^2,
 \end{equation}
 where $\kappa$ is the curvature of the curve $x\mapsto\varphi(x)$.

$S_1$ generalizes the real part of the analytic Schwarzian, while the connection we need between $S_1$ and the Schwarzian for harmonic maps is
\[
S_1\widetilde{f}(x) \le \Re\{\Sch f(x)\} + e^{2\sg(x)}|K(\widetilde{f}(x))|, \quad -1 < x <1;
\]
see Lemma 1 in \cite{cdo:injective-lift}. Thus if $f$ satisfies \eqref{eq:nehari-harmonic} then
\begin{equation} \label{eq:S1-bound}
S_1\widetilde{f}(x) \le \frac{2}{(1-x^2)^2}, \quad -1 < x < 1.
\end{equation}

 With all this as background, our first lemma is fairly straightforward.

 \begin{lemma} \label{lemma:h-convex}
 Let $f$ satisfy \eqref{eq:nehari-harmonic}, with lift $\widetilde{f}$, and let $T$ be a M\"obius transformation of $\overline{\mathbb{R}^3}$. The function
 \[
 u_{T \circ \widetilde{f}}(z) = \frac{1}{\sqrt{(1-|z|^2)e^{\tau(z)}}}, \quad e^\tau=(\|T'\|\circ \widetilde{f})e^\sg,
 \]
 is hyperbolically convex in $\D$.
 \end{lemma}
While $(T\circ \widetilde{f})(\D) = T(\Sigma)$ is generally not a minimal surface, $T$ \emph{is} a conformal mapping of $\overline{\mathbb{R}^3}$ and $e^{2\tau}|dz|^2$ is the corresponding conformal metric on $\D$.

\begin{proof}[Proof of Lemma \ref{lemma:h-convex}]
 Since $u_{\widetilde{f}}$ is hyperbolically convex and we can consider $\widetilde{f} \circ M$ for any M\"obius transformation of $\D$ onto itself, it suffices to show that $u_{T\circ\widetilde{f}}$ is hyperbolically convex along the diameter $-1 <x<1$. The argument proceeds by comparing coefficients in two second-order differential equations.

 Let $\varphi(x) = (T \circ \widetilde{f})(x)$.  From M\"obius invariance and  \eqref{eq:S1-bound},
 \[
 S_1\varphi(x) = S_1 \widetilde{f}(x) \le \frac{2}{(1-x^2)^2}, \quad -1 < x < 1.
 \]
Now with $v(x) = |\varphi'(x)|= e^{\tau(x)}$, as above, from \eqref{eq:S1-curvature}
\begin{equation} \label{eq:v-and-S1}
\left(\frac{v'(x)}{v(x)}\right)'  -\frac{1}{2}\left(\frac{v'(x)}{v(x)}\right)^2 \le S_1\varphi(x)  \le \frac{2}{(1-x^2)^2}. %
\end{equation}
Let $2p$ denote the left-hand side, so that
\begin{equation} \label{eq:p<=}
2p(x) \le \frac{2}{(1-x^2)^2}, \quad -1 < x <1.
\end{equation}
The function $V=v^{-1/2}$ satisfies the differential equation
\begin{equation} \label{eq:V''}
V''+pV = 0
\end{equation}
and the function
\begin{equation} \label{eq:W}
W(x) = \frac{V(x)}{\sqrt{1-x^2}}
\end{equation}
is precisely $u_{T\circ \widetilde{f}}$ restricted to $-1 <x <1$. If we give $-1 <x <1$ its hyperbolic parametrization,
\[
s=\frac{1}{2}\log\frac{1+x}{1-x}, \quad x(s) = \frac{e^{2s}-1}{e^{2s}+1}, \quad x'(s) = 1-x(s)^2,
\]
 a calculation produces
\[
\frac{d^2}{ds^2}W = \left(\frac{1}{(1-x^2)^2}-p(x)\right)(1-x^2)^2W(x), \quad x=x(s),
\]
and appealing to \eqref{eq:p<=} shows this is nonnegative.
\end{proof}


The topological condition that  $\widetilde{f}$ be injective on $\partial\D$ has an analytical consequence on critical points that is important for much of our work.

\begin{lemma} \label{lemma:critical-point}
If $f$ satisfies \eqref{eq:nehari-harmonic} and is injective on $\partial\D$ then function $u_{T \circ \widetilde{f}}$ has at most one critical point in $\D$.
\end{lemma}

\begin{proof}
Suppose that $u_{T\circ\widetilde{f}}$ has two critical points. Composing $\widetilde{f}$ with a M\"obius transformation of $\D$ onto itself we may locate the critical points at $0$ and $a$, $0<a<1$. By convexity these must give absolute minima of $u_{T\circ\widetilde{f}}$ in $\D$, and the same must be true of $u_{T\circ\widetilde{f}}(x)$ for $0\le x \le a$. Hence $u_{T\circ \widetilde{f}}$ is constant on $[0,a]$ and thus constant on $(-1,1)$ because it is real analytic there.

It follows that the function $v(x) = e^{\tau(x)}$ is a constant multiple of $1/(1-x^2)^2$. But then $V(x)= v(x)^{-1/2}$ is constant multiple of $\sqrt{1-x^2}$, and from the differential equation \eqref{eq:V''} we conclude that $p(x) = 1/(1-x^2)^2$. In turn, from \eqref{eq:S1-curvature} and  \eqref{eq:v-and-S1} this forces the curvature $\kappa$ to vanish identically. Thus $T\circ \widetilde{f}$ maps the interval $(-1,1)$ onto a line with speed $\|\varphi'(x)\|=v(x) = 1/(1-x^2)$, and so $\varphi(1) = \varphi(-1) = \infty$. This  violates the assumption that $\widetilde{f}$, hence $T\circ \widetilde{f}$, is injective on $\partial\D$.
\end{proof}

Continuing with the same assumptions, we now show what happens when there is exactly one critical point.
\begin{lemma} \label{lemma:bounded-equivalences}
Let $f$ satisfy \eqref{eq:nehari-harmonic} and be injective on $\partial\D$. Let $T$ be a M\"obius transformation of $\overline{\mathbb{R}^3}$. The following are equivalent:
\begin{list}{}{\setlength\leftmargin{.25in}}
\item[$(i)$]  $u_{T\circ\widetilde{f}}$ has a critical point.
\item[$(ii)$]  $(T\circ \widetilde{f})(\D)$ is bounded.
\item[$(iii)$] $u_{T\circ \widetilde{f}}(re^{i\theta})$ is eventually increasing along each radius $[0,e^{i\theta})$.
\item[$(iv)$]  $u_{T\circ\widetilde{f}}(z) \rightarrow \infty$ as $|z|\rightarrow 1$.
\end{list}
\end{lemma}

In the proof of this lemma, and elsewhere, we will have occasion to use M\"obius inversions. Following Ahlfors we write
\[
J(x) = \frac{x}{\|x\|^2}
\]
and for the derivative
\begin{equation} \label{eq:DJ}
J'(x) = \frac{1}{\|x\|^4}(\|x\|^2{\rm Id}-2Q(x)),
\end{equation}
where
\[
Q(x)_{ij} = {x_ix_j}.
\]
and Id is the identity. From this and $Q(x)^2 = \|x\|^2Q(x)$ one has
\begin{equation} \label{eq:norm-DJ}
\|J'(x)\|= \frac{1}{\|x\|^2}.
\end{equation}

\begin{proof}[Proof of Lemma \ref{lemma:bounded-equivalences}]
If $(iv)$ holds there is an interior minimum so $(iv) \implies (i)$ is immediate.

Suppose $(i)$ holds. We may assume the critical point is at the origin. The value $u_{T\circ\widetilde{f}}(0)$ is the absolute minimum for  $u_{T\circ\widetilde{f}}$ in $\D$ and so
\[
e^{\tau(z)} \le \frac{e^{\tau(0)}}{1-|z|^2}, \quad z \in \D.
\]
Thus  $\tau$ remains finite in $\D$ and $\infty$ cannot be a point on $T(\Sigma)$.

To show that $T(\Sigma)$ is bounded we first work along $[0,1)$. The hyperbolically convex function $W(x) =u_{T\circ\widetilde{f}}(x)$ in \eqref{eq:W} cannot be constant because $0$ is the unique critical point. Hence if $x(s)$ is the  hyperbolic arclength parametrization of $[0,1)$ with $x(0)=0$ then
\[
\frac{d}{ds} W(x(s)) \ge a, \quad W(x(s)) \ge a s + b,
\]
for some $a, b >0$ and all $s \ge s_0 >0$. From this
\[
\begin{aligned}
v(x) = \frac{1}{V(x)^2} & \le \frac{1}{(1-x^2)\left(\frac{a}{2}\log\frac{1+x}{1-x}+b\right)^2}\\
&= -\frac{1}{a}\frac{d}{dx}\left(\frac{1}{\frac{a}{2}\log\frac{1+x}{1-x}+b}\right).
\end{aligned}
\]
Therefore
\[
\int_0^1 e^{\tau(x)}\,dx = \int_0^1v(x)\,dx <\infty,
\]
with a bound depending only on $a$, $b$, $s_0$, and $(T\circ \widetilde{f})(1)$ is finite.

This argument can be applied on every radius $[0,e^{i\theta})$, and by compactness the corresponding numbers $a_\theta, b_\theta, s_\theta$ can be chosen positive independent of $\theta$. This proves that $T\circ\widetilde{f}$ is bounded, and hence that $(i) \implies (ii)$.

For $(ii)$ $\implies$ ${(iii)}$ we can first rotate and assume $e^{i\theta} = 1$. In the notation above, we need to show for some $x_0 >0$ that $W(x)$ is increasing for $x_0\le x <1$. 

We have to follow $T$ by an inversion, so to simplify the notation let $\widetilde{f}_1= T\circ \widetilde{f}$ and
$u_{\widetilde{f}_1}(z) = ((1-|z|^2)e^{\tau(z)})^{-1/2}$. For $w_0$ to be determined let
\[
I(w) = \frac{w-w_0}{\|w-w_0\|^2},
\]
and write
\[
\widetilde{f_2}=I\circ\widetilde{f_1},\quad
u_{\widetilde{f_2}}(z) = \frac{1}{\sqrt{(1-|z|^2)e^{\nu(z)}}}, \quad \text{and} \quad W_2(x) = u_{\widetilde{f_2}}(x),\, x \in (-1,1).
\]
Again, we know that $W_2(x(s))$ is convex, where $s$ is the hyperbolic arclength parameter.

From \eqref{eq:norm-DJ},
\[
e^{\nu(z)}=\frac{e^{\tau(z)}}{\|\widetilde{f_1}(z) - w_0\|^2} \quad \text{or} \quad \nu(z) = \tau(z) - \log\|\widetilde{f_1}(z)-w_0\|^2,
\]
and therefore
\begin{equation} \label{eq:grad-equation}
\nabla\nu(0) = \nabla\tau(0) + \frac{2}{\|w_0\|^2}\left(\left\langle \frac{\partial \widetilde{f}_2}{\partial x}(0), w_0\right\rangle, \left\langle \frac{\partial \widetilde{f_2}}{\partial y}(0),w_0\right\rangle\right).
\end{equation}
But also 
\[
\nabla u_{\widetilde{f}_2}(0) = -\frac{1}{2}\nabla \nu(0),
\]
and from this equation and \eqref{eq:grad-equation} it is clear we can choose $w_0$ to make
\[
W_2'(0) = a>0.
\]
Convexity then ensures $W_2(x(s)) \ge a s$.

To work back to $W$, write
\begin{equation} \label{eq:f1-and-f2}
\widetilde{f_1} = \frac{\widetilde{f_2}}{\|\widetilde{f_2}\|^2}+w_0,
\end{equation}
whence
\[
\|D\widetilde{f_1}\| =\frac{\|D\widetilde{f_2}\|}{\|\widetilde{f_2}\|^2},
\]
and
\[
W=W_2\|\widetilde{f}_2\|.
\]
The assumption we make in $(ii)$ is that $\widetilde{f_1}(\D) = (T\circ \widetilde{f})(\D)$ is bounded, and \eqref{eq:f1-and-f2} thus implies that  $\|\widetilde{f}_2\|\ge \delta>0$. Therefore $W(x(s)) \ge a\delta s$. By convexity, there is an $x_0>0$ so that $W(x)$ is increasing for $x_0 \le x <1$.  This completes the proof that $(ii) \implies (iii)$.

Finally, if $(iii)$ holds then for each $\theta$ there exists $0 < r_\theta <1$ such that
\[
\frac{\partial}{\partial r} u_{T\circ \widetilde{f}}(r_\theta e^{i\theta}) \ge a_\theta >0.
\]
 By compactness the $r_\theta$ can be chosen bounded away from $1$  and the $a_\theta$ bounded away from $0$. By hyperbolic convexity, along the tail of each radius $u_{T \circ\widetilde{f}}(r(s)e^{i\theta})$ is uniformly bounded below by a linear function of the hyperbolic arclength parameter $s$, which tends to $\infty$ as $r= r(s) \rightarrow 1$.
\end{proof}

We conclude this section with some remarks on introducing the function $u_{\widetilde{f}}$. Let
\[
\lambda_\D(z)^2|dz|^2= \frac{1}{(1-|z|^2)^2}|dz|^2
\]
be the Poincar\'e metric for $\D$ (curvature $-4$) and let  $\lambda_\Sigma^2\,\mathbf{g}_0$ be the conformal metric on $\Sigma$ with
\[
\widetilde{f}^*(\lambda_\Sigma^2\,\mathbf{g}_0) = \lambda_\D^2 |dz|,
\]
so that $\widetilde{f}$ is an isometry. 
Since $\widetilde{f}^*(\mathbf{g}_0) = e^{2\sg}|dz|^2$ we have
\begin{equation} \label{eq:lambda-Sigma}
(\lambda_\Sigma \circ \widetilde{f})(z)  = \frac{1}{(1-|z|^2)e^{\sg(z)}} \quad \text{or} \quad \lambda_\Sigma \circ \widetilde{f} = e^{- \sigma}\lambda_\D,
\end{equation}
and
\[
u_{\widetilde{f}} = (\lambda_\Sigma\circ \widetilde{f})^{1/2}.
\]

\smallskip

If $f$ is analytic and injective in $\D$ and the plane domain $\Omega =f(\D)$ replaces $\Sigma$, then $\lambda_\Sigma = \lambda_\Omega$ is the Poincar\'e metric on $\Omega$ and $u_f = (\lambda_\Omega \circ f)^{1/2}$. In \cite{cop:john-nehari} it was shown that the hyperbolic convexity of $\lambda_{T(\Omega)}^{1/2}$ for any M\"obius transformation $T$ is a characteristic property of functions satisfying the Nehari condition \eqref{eq:nehari-2}. Lemma \ref{lemma:h-convex} is an analog of this for harmonic maps. We will use $\lambda_\Sigma$ to write various identities, inequalities, etc., in forms intrinsic to $\Sigma$.


\section{Circles and Reflections} \label{section:circles-reflections}

We continue to assume that $f$ satisfies the injectivity condition \eqref{eq:nehari-harmonic} with a lift $\widetilde{f}$ mapping $\D$ to the minimal surface $\Sigma\subset \mathbb{R}^3$, and also that $\widetilde{f}$ is injective on $\partial\D$. The purpose of this section is to define a reflection $R\colon \Sigma \longrightarrow \overline{\mathbb{R}^3}\setminus \overline{\Sigma}$ that provides a continuous, injective extension of $\widetilde{f}$.
To extend $\Sigma$ beyond its boundary we use a family of Euclidean circles (possibly including a line) each orthogonal to $\Sigma$. They are defined by the following lemma, which depends on properties of the function $u_{T\circ \widetilde{f}}$ established in the preceding section.

\begin{lemma} \label{lemma:circles}
For each $w \in \Sigma$ there is a unique Euclidean circle $C_w \subset \overline{\mathbb{R}^3}$ with the following properties:
\begin{list}{}{\setlength\leftmargin{.25in}}
\item[$(i)$] $C_{w}$ is orthogonal to $\Sigma$ at $w$;
\item[$(ii)$] $C_{w} \cap \Sigma = \{w\}$;
\item[$(iii)$] Let $z_0\in \D$ and $w_0=f(z_0)$. A point $w_1$ lies on $C_{w_0} \setminus \{w_0\}$ if and only if  $u_{I \circ \widetilde{f}}$ has a critical point at $z_0$, where
\[
I(w) = \frac{w-w_1}{\|w-w_1\|^2},\quad w \in \overline{\mathbb{R}^3}.
\]
\end{list}
If $u_{\tilde{f}}$ has a critical point at $z_0$ then $C_{w_0}$ is a line satisfying $(i)$ and $(ii)$.
\end{lemma}
Briefly, when referring to part $(iii)$ we say that inversion about any point in $C_{w_0}$ other than $w_0$ produces a critical point for $u_{I\circ \widetilde{f}}$ at  $z_0=\widetilde{f}^{-1}(w_0) \in \D$. Observe that if $I\circ \widetilde{f}$ produces a critical point for $u_{I \circ \widetilde{f}}$ at $z_0$ then so does any further \emph{affine} change $A \circ  I \circ \widetilde{f}$. We will need this later.

\begin{proof}
We begin by determining the conditions under which $u_{I\circ \widetilde{f}}$ has a critical point when $I$ is an inversion.  This recapitulates some of the calculations in the proof of Lemma \ref{lemma:bounded-equivalences}.

Consider first the case $z_0=0$. We can also assume that $\widetilde{f}(0)=0$, and we let $T_0\Sigma$ denote the tangent plane to $\Sigma$ at $0$. Computing from the definition of $u_{T\circ\widetilde{f}}$ we have, as in \eqref{eq:grad-equation}, that  $\nabla u_{T\circ\widetilde{f}}(0) = 0$ if and only if $\nabla \tau(0)=0$, and this is for any M\"obius transformation $T$.   Specializing to the inversion
\begin{equation} \label{eq:w_0-inversion}
I(w) = \frac{w-w_1}{\|w-w_1\|^2},\quad       (I\circ \widetilde{f})(z) = \frac{\widetilde{f}(z)-w_1}{\|\widetilde{f}(z)-w_1\|^2}, \quad w_1 \ne 0,
\end{equation}
gives for $u_{I\circ \widetilde{f}}$ that
\[
e^{\tau(z)} = \frac{e^{\sg(z)}}{\|\widetilde{f}(z)-w_1\|^2} \quad \text{or} \quad \tau(z) = \sg(z) - \log\|\widetilde{f}(z)-w_1\|^2.
\]
Thus
\begin{equation} \label{eq:grad-tau}
\nabla\tau(0) = \nabla\sg(0) + \frac{2}{\|w_1\|^2}(\langle\widetilde{f}_x(0), w_1\rangle\,,\, \langle \widetilde{f}_y(0),w_1\rangle),
\end{equation}
and $\nabla\tau(0) = 0$ when
\begin{equation} \label{eq:w_0-orthogonal}
\frac{1}{\|w_1\|^2}\langle\widetilde{f}_x(0),w_1\rangle = -\frac{1}{2} \sg_x(0) \quad \text{and}\quad \frac{1}{\|w_1\|^2}\langle\widetilde{f}_y(0),w_1\rangle = -\frac{1}{2} \sg_y(0).
\end{equation}
Since $\widetilde{f}$ is conformal
\[
\langle \widetilde{f}_x(0) , \widetilde{f}_y(0) \rangle =0 \quad \text{and} \quad \|\widetilde{f}_x(0)\|^2=\|\widetilde{f}_y(0)\|^2 = e^{2\sg(0)}.
\]
Then \eqref{eq:w_0-orthogonal} says exactly that the point $w_1/\|w_1\|^2$ lies on a  line orthogonal to $T_0\Sigma$ through the point
\[
\zeta= -\frac{1}{2}e^{-2\sg(0)}\left\{\sg_x(0)\widetilde{f}_x(0)+\sg_y(0)\widetilde{f}_y(0)\right\}.
\]
on $T_0\Sigma$. Call this line $L_0$; it depends only on the various data at $0$.

The inversion
$
J(w) = {w}/{\|w\|^2}
$
leaves the tangent plane $T_0\Sigma$ invariant and interchanges $0$ and $\infty$, where $L_0$ and $T_0\Sigma$ meet a second time orthogonally. That is, if we put
$
C_{0} = J(L_0)
$
then $w_1\in C_0$ and $C_0$ is orthogonal to $\Sigma$ at $0$ and also, generically, orthogonal to $T_0\Sigma$ at some other finite point (which we will determine). The exceptional case is when $L\cap T_0\Sigma = \{0\}$, which occurs when $\nabla\sg(0)=0$. Then $\nabla\tau(0)=0$ and $u_{\widetilde{f}}$ already has a critical point at $0$. In this case $C_0=J(L_0) = L_0$.  This proves parts $(i)$ and $(iii)$ of the lemma for $z_0=0$.

Part (ii) of the lemma, for $z_0=0$, follows from Lemma \ref{lemma:bounded-equivalences}. Indeed, if the inversion \eqref{eq:w_0-inversion} produces a critical point for $u_{I\circ \widetilde{f}}$ (at $0$) then $I(\Sigma)$ is bounded, and hence $w_1$ cannot lie on $\Sigma$.

Finally, to pass from $0$ to an arbitrary point $z_0 \in \D$,  consider
\[
\widetilde{f}_1(z) = f\left(\frac{z+z_0}{1+\overline{z_0}z}\right) - \widetilde{f}(z_0).
\]
By Schwarz' lemma
\[
u_{\widetilde{f}_1(z)} = u_{\widetilde{f}}\left(\frac{z+z_0}{1+\overline{z_0}z}\right),
\]
hence $u_{I\circ\widetilde{f}_1}$ has a critical point at $0$ if and only if $u_{I\circ\widetilde{f}}$ has a critical point at $z_0$. The statements $(i)$, $(ii)$ and $(iii)$ then follow from the previous analysis.
\end{proof}

By means of this construction, each point $w$ on $\Sigma$ is associated to a point $w^*$ outside $\Sigma$  on the tangent plane $T_w\Sigma$, namely the other point where $C_w$ meets $T_w\Sigma$. The points $w$ and $w^*$ are endpoints of the diameter of $C_w$ that lies in $T_w\Sigma$. We write $w^*=R(w)$, or $R_\Sigma(w)$,  and refer to $w^*$ as the reflection of $w$. In Section \ref{section:extension} we will show that $R$ fixes $\partial\Sigma$ pointwise.
Note also that the arguments used to define the reflection of $\Sigma$ can be applied to define  the reflection $R_{\Sigma'}$ of any surface $\Sigma'=T(\Sigma)$, $T$ a M\"obius transformation, using the function $u_{T\circ\widetilde{f}}$. It is not true, however, that $R_{\Sigma'}\circ T = T\circ R_\Sigma$. We will return to this at the end of this section.

It is a consequence of Lemma \ref{lemma:critical-point} that $R$ is injective.

\begin{lemma} \label{lemma:injective}
If $w \ne w'$ then $C_{w}\cap C_{w'} = \emptyset$. Hence $R$ is injective.
\end{lemma}

\begin{proof}
The circles meet $\Sigma$ only at the distinct points $w$ and $w'$.  If there is a point $w_1\in C_{w}\cap C_{w'}$ it is not on $\Sigma$ and  the inversion $I(w) =(w-w_1)/(\|w-w_1\|^2)$ produces critical points for $u_{I\circ \widetilde{f}}$ at distinct points $z=\widetilde{f}^{-1}(w)$ and $z'=\widetilde{f}^{-1}(w')$ in $\D$. This is impossible by Lemma \ref{lemma:critical-point}.
\end{proof}

 It is not difficult to find a formula for $w^*=R(w)$. The vectors
 \[
 X(z)= e^{-\sg(z)}\widetilde{f}_x(z) ,\quad Y(z)=e^{-\sg(z)}\widetilde{f}_y(z)
 \]
 are an orthonormal basis for $T_w\Sigma$, $w=\widetilde{f}(z)$. 
 Again, first take $z=0$ and $w=\widetilde{f}(z)=0$.  Since $w^* \in C_0$ the equations  \eqref{eq:w_0-orthogonal} apply to $w^*$ and from these
\[
\frac{1}{\|w^*\|^2}\langle w^*,X\rangle = -\frac{1}{2}e^{-\sg(0)}\sg_x(0), \quad \frac{1}{\|w^*\|^2}\langle w^*,Y\rangle = -\frac{1}{2}e^{-\sg(0)}\sg_y(0).
\]
This leads easily to
\[
w^* = -\frac{2e^{\sg(0)}}{\|\nabla\sg(0)\|^2}\left\{\sg_x(0)X(0) + \sg_y(0)Y(0)\right\}.
\]

The formula when $z$ is any point in $\D$ and $w=\widetilde{f}(z)$ is obtained by renormalizing $\widetilde{f}$ as in the proof of the lemma, including a translation by $\widetilde{f}(z)$. The result is
\begin{equation} \label{eq:w^*}
w^* = w + \frac{e^{\sg(z)}\alpha(z)}{\alpha(z)^2+\beta(z)^2}X(z) + \frac{e^{\sg(z)}\beta(z)}{\alpha(z)^2+\beta(z)^2}Y(z)
\end{equation}
where
\begin{equation} \label{eq:A,B}
\alpha(z) = \frac{x}{1-|z|^2}-\frac{1}{2}\sg_x(z), \quad \beta(z) = \frac{y}{1-|z|^2}-\frac{1}{2}\sg_y(z), \quad z=x+iy.
\end{equation}
One can also verify
\[
\nabla\log u_{\widetilde{f}}(z) = (\alpha(z),\beta(z)).
\]
The function $u_{\widetilde{f}}$ has a critical point at $z$ precisely when $\alpha(z)=\beta(z)=0$, in which case $w^*=\infty$. Note as well that the diameter of $C_w$ is
\begin{equation} \label{eq:diameter-Cw}
\|w^*-w\|= \frac{e^{\sg(z)}}{\|\nabla \log u_{\widetilde{f}}(z)\|}.
\end{equation}

\smallskip

Furthermore, we can write the reflection in a form intrinsic to the surface.
\begin{lemma} \label{lemma:R}
The reflection $w^*=R(w)$ is given by
\begin{equation} \label{eq:R-gradient}
R(w) = w+2J(\nabla\log\lambda_\Sigma(w)),
\end{equation}
where $J(w) = w/\|w\|^2$.
\end{lemma}

\begin{proof}
Recall from \eqref{eq:lambda-Sigma} the conformal metric $\lambda_\Sigma^2\,\mathbf{g}_0$ on $\Sigma$ that is isometric to the Poincar\'e metric $|dz|^2/(1-|z|^2)^2$ on $\D$ and the relation
\[
\log(\lambda_\Sigma \circ \widetilde{f})(z) =  -\log(1-|z|^2) -\sg(z).
\]
Then with \eqref{eq:A,B},
\[
\begin{aligned}
&\frac{1}{2}\frac{\partial}{\partial x}\log(\lambda_\Sigma \circ \widetilde{f})  =  \frac{x}{1-|z|^2}-\frac{1}{2}\sg_x =\alpha\\
& \frac{1}{2}\frac{\partial}{\partial y}\log(\lambda_\Sigma \circ \widetilde{f})  =  \frac{y}{1-|z|^2}-\frac{1}{2}\sg_y=\beta .
\end{aligned}
 \]
Now let $\nabla \log\lambda_\Sigma$ be the gradient with respect to the Euclidean metric on $\Sigma$. As  a vector field on $\Sigma$ we can write, with $w=\widetilde{f}(z)$,
\[
\begin{aligned}
\nabla\log\lambda_\Sigma(w) &= e^{-\sg(z)}\left\{\frac{\partial}{\partial x} (\log\lambda_\Sigma \circ \widetilde{f})(z)X(z)+\frac{\partial}{\partial y} (\log\lambda_\Sigma \circ \widetilde{f})(z)Y(z)\right\}\\
&= 2e^{-\sg(z)}\{\alpha(z)X(z)+\beta(z)Y(z)\}
\end{aligned}
\]
and
\[
\|\nabla\log\lambda_\Sigma(w)\|^2 = 4e^{-2\sg(z)}(\alpha(z)^2+\beta(z)^2).
\]
Using  the inversion $J(w) = w/\|w\|^2$ we thus have
\[
2J(\nabla\log\lambda_\Sigma(w)) =  \frac{e^{\sg(z)}\alpha(z)}{\alpha(z)^2+\beta(z)^2}X(z) + \frac{e^{\sg(z)}\beta(z)}{\alpha(z)^2+\beta(z)^2}Y(z)
\]
and
\[ 
R(w) = w+2J(\nabla\log\lambda_\Sigma(w)).
\]
as stated.
\end{proof}

Finally we consider a conformal invariance property of the construction. This will be important in the next section when we show that $\widetilde{f}$ and its extension match on $\partial\D$.

One cannot expect $R_\Sigma$ to be conformally natural, meaning that
\[
R_{\Sigma'}(T(w)) = T(R_\Sigma(w))
\]
for a M\"obius transformation $T$ with $T(\Sigma) = \Sigma'$, since $w^*=R_\Sigma(w)$ is defined at each point  $w$ via the tangent {plane} to the surface and under a M\"obius transformation this plane may become a sphere.  
However, 
if the families of circles $\{C_w \colon w\in \Sigma\}$ and $\{C_\omega \colon \omega \in \Sigma'\}$ define the reflections for the surfaces $\Sigma$ and $\Sigma'$, respectively, then $T(C_w) = C_{T(w)}$.
We can describe this degree of conformal invariance succinctly by introducing
\[
\mathcal{C}_\Sigma = \bigcup_{w\in \Sigma} C_w.
\]
Then
\begin{equation} \label{eq:conformal-invariance}
T(\mathcal{C}_\Sigma) = \mathcal{C}_{T(\Sigma)}.
\end{equation}

To show this,  note that as $w=\widetilde{f}(z)$ varies over $\Sigma$, the circles $T(C_w)$ clearly have properties $(i)$ and $(ii)$ of          Lemma~\ref{lemma:circles} 
for the surface $\Sigma'$. Take a point $w_0=\widetilde{f}(z_0)$, determining the circle $C_{w_0}$, and let $\omega_0=T(w_0)$.  The question is whether for any $\omega_1 =T(w_1) \in T(C_{w_0}) \setminus\{\omega_0\}$ the inversion
\[
I(\zeta) = \frac{\zeta - \omega_1}{\|\zeta - \omega_1\|^2}
\]
produces a critical point for $u_{I\circ T\circ \widetilde{f}}$ at $z_0$.
But  the map
\[
(I \circ T)(v) = \frac{T(v) - \omega_1}{\|T(v) - \omega_1\|^2}
\]
is a M\"obius transformation sending $w_1$ to $\infty$, as is the inversion
\[
K(v) = \frac{v-w_1}{\|v-w_1\|^2}.
\]
It follows that
\[
(I \circ T )(v) = 
 (A \circ K)(v) 
\]
for an affine transformation $A$. Now the circle $C_{w_0}$ has the property  that the inversion $K(v) = (v-w_1)/\|v-w_1\|^2$ produces a critical point for $u_{K\circ f}$ at $z_0$, and since $A$ is affine, $I\circ T=A \circ K$ produces a critical point for $u_{I\circ T\circ \widetilde{f}}$ at $z_0$ as we were required to show.  We conclude that the circles $T(C_w)$ for the surface $\Sigma'$ have the properties of the circles in Lemma \ref{lemma:circles}, and  that $T(C_w)= C_{T(w)}$. 

In addition to the conformal invariance expressed by  \eqref{eq:conformal-invariance} we have
\begin{equation} \label{eq:C-fills-space}
\mathcal{C}_\Sigma \cup \partial\Sigma = \overline{\mathbb{R}^3},
\end{equation}
and by Lemma \ref{lemma:injective} this is a disjoint union. We will not need \eqref{eq:C-fills-space} but we consider it an important feature of the construction. To prove it, observe that if $w_1\not\in \overline{\Sigma}=\Sigma \cup \partial\Sigma$ then the inversion
\[
I(w) = \frac{w-w_1}{\|w-w_1\|^2}
\]
has the property that $I(\Sigma)$ is bounded. It follows by Lemma~\ref{lemma:bounded-equivalences} that $u_{I\circ\widetilde{f}}$ has a critical point, and  by Lemma~\ref{lemma:circles} that $w_1$ lies on some $C_{w_0}\setminus\{w_0\}$.

\section{Definition of the Extension and Proof of Theorem \ref{theorem:extension}, Part (a)} \label{section:extension}

With assumptions and notations as before, we define
\begin{equation} \label{eq:F}
\widetilde{F}(z)=
\begin{cases}
\widetilde{f}(z), & \quad z \in\overline{\D},\\
R({\widetilde{f}}(1/\bar{z})), & \quad z \in \overline{\mathbb{C}}\setminus \overline{\D}.
\end{cases}
\end{equation}
To prove that $\widetilde{F}$ defines an extension of $\widetilde{f}$ we must show that $\widetilde{f}$ and $R\circ \widetilde{f}$ match continuously along $\partial\D$.

\begin{lemma} \label{lemma:boundary-values}
Let $z \in \D$ and let $d$ denote the spherical metric on $\overline{\mathbb{R}^3}$. Then
\[
d(\widetilde{f}(z), R(\widetilde{f}(z))) \rightarrow 0, \quad |z| \rightarrow 1.
\]
\end{lemma}

\begin{proof}
 We divide the proof into the cases when $u_{\widetilde{f}}$ has one critical point and when it has none. We work in the spherical metric because, first, $\widetilde{f}$ has a spherically continuous extension, and second, when $u_{\widetilde{f}}$ has no critical points we have to allow for shifting $\widetilde{f}$ by a M\"obius transformation.

Suppose $u_{\widetilde{f}}$ has a unique critical point, which we can take to be at $0$. The proof of Lemma~\ref{lemma:bounded-equivalences} shows that there is an $a >0$ such that along any radius $[0, e^{i\theta})$
\[
(1-r^2)\frac{\partial}{\partial r} u_{\widetilde{f}}(re^{i\theta}) \ge a
\]
for all $r \ge r_0>0$. (This corresponds to $dW/ds \ge a$ in the proof of Lemma~\ref{lemma:bounded-equivalences}, where $s$ is the hyperbolic arclength parameter.) From this it follows that
\[
(1-|z|^2)\|\nabla u_{\widetilde{f}}(z)\| \ge a >0,
\]
for all $|z| \ge r_0 >0$.

From \eqref{eq:diameter-Cw}
\[
\begin{aligned}
\|R(\widetilde{f}(z)) - \widetilde{f}(z)\| &= \frac{e^{\sg(z)}}{\| \nabla \log u_{\widetilde{f}}(z)\|}
=\frac{u_{\widetilde{f}}(z)e^{\sg(z)}}{\|\nabla u_{\widetilde{f}}(z)\|}\\
&= \frac{1}{u_{\widetilde{f}}(z)}\frac{1}{(1-|z|^2)\|\nabla u_{\widetilde{f}}(z)\|}.
\end{aligned}
\]
This tends to $0$ as $|z|\rightarrow 1$ because $u_{\widetilde{f}}$ becomes infinite (Lemma \ref{lemma:bounded-equivalences}) and $(1-|z|^2)\|\nabla u_{\widetilde{f}}(z)\|$ stays bounded below. Geometrically,  the diameter of  $C_{\widetilde{f}(z)}$ tends to $0$ as $|z|$ increases to $1$.

Next, supposing that $u_{\widetilde{f}}$ has no critical point we produce one. That is, let $T$ be a M\"obius  transformation so that $u_{T \circ \widetilde{f}}$ has a critical point at $0$. The preceding argument can be repeated verbatim to conclude that
\begin{equation} \label{eq:diameter-shrinks}
\|R(T(\widetilde{f}(z))) - T(\widetilde{f}(z))\|\rightarrow 0 \quad \text{as} \quad  |z| \rightarrow 1.
\end{equation}
 If $R$ were conformally natural, if we knew that $R \circ T= T \circ R$, then  we would be done.  Instead, we argue as follows.

Let $z\in \D$, $z \ne 0$. The length $\|R(T(\widetilde{f}(z))) - T(\widetilde{f}(z))\|$ is the diameter of the circle $C_{T(\widetilde{f}(z))}$ based at $T(\widetilde{f}(z))$  that defines the reflection for the surface $T(\Sigma)$, and it tends to $0$ by \eqref{eq:diameter-shrinks}. But now, if $C_{\widetilde{f}(z)}$ is the circle based at $\widetilde{f}(z)$, for the surface $\Sigma$,  then the reflected point $R(\widetilde{f}(z))$ is also on this circle (diametrically opposite $\widetilde{f}(z)$) and then $T(R(\widetilde{f}(z))) \in C_{T(\widetilde{f}(z))}$. Therefore $\|T(R({\widetilde{f}}(z))) - T(\widetilde{f}(z))\|\rightarrow 0$ as $|z|\rightarrow 1$, whence in the spherical metric $d(R({\widetilde{f}}(z)), \widetilde{f}(z))$ 
tends to $0$ as well and the proof is complete.
\end{proof}

Combining Lemmas \ref{lemma:injective} and \ref{lemma:boundary-values} proves part (a) of Theorem \ref{theorem:extension}; the mapping $\widetilde{F}$ defined in \eqref{eq:F} is a continuous injective extension of $\widetilde{f}$. Furthermore, the formulas make clear that $\widetilde{F}$ is real-analytic off $\partial\D$.

\subsection*{Remarks} When $f$ is analytic in $\D$ the Ahlfors-Weill extension extension can be written as
\[
F(z) =
\begin{cases}
f(z), \quad z \in \overline{\D},\\
f(\zeta) + \displaystyle{\frac{(1-|\zeta|^2)f'(z)}{\bar{\zeta} - \frac{1}{2}(1-|\zeta|^2)\frac{f''(\zeta)}{f'(\zeta)}}}, \quad \zeta = 1/\bar{z}, \,z \in \mathbb{C} \setminus \overline{\D}.
\end{cases}
\]
Ahlfors and Weill did not express it in this form; see \cite {co:aw}. Alternatively, if $\lambda_\Omega|dw|$ is the Poincar\'e metric on $\Omega = f(\D)$ then
\[
F(z) =
\begin{cases}
f(z), \quad z \in \overline{\D},\\
f(\zeta) +\displaystyle{\frac{1}{\partial_w\log\lambda_\Omega(f(\zeta))}}, \quad \zeta = 1/\bar{z}, \,z \in \mathbb{C} \setminus \overline{\D}.
\end{cases}
\]
The equation \eqref{eq:R-gradient} for the reflection gives exactly
\begin{equation} \label{eq:aw-extension-gradient}
R(w) = w + \frac{1}{\partial_w\log\lambda_\Omega(w)}
\end{equation}
 when $f$ is analytic.

The Ahlfors-Weill reflection is conformally natural: If $T$ is a M\"obius transformation of $\overline{\mathbb{C}}$ and $T(\Omega) = \Omega'$ then
\[
R_{\Omega'} \circ T = T \circ R_\Omega.
\]
From the perspective of the present paper this is because all tangent planes $T_z(\Omega)$ to $\Omega$ can be identified with $\mathbb{C}$, which is preserved by the extensions to $\overline{\mathbb{R}^3}$ of the M\"obius transformations.

 The  reflection defining the Ahlfors-Weill extension was expressed in a form like \eqref{eq:aw-extension-gradient} also by Epstein \cite{epstein:reflections} in his penetrating geometric study of Nehari's and related theorems. Still another interesting geometric construction, using Euclidean circles of curvature, was given by Minda \cite{minda:reflections}.

\section{Quasiconformality of the Reflection and Proof of Theorem \ref{theorem:extension}, Part (b)} \label{section:qc-calculation}

We now assume that $f$ satisfies
\begin{equation} \label{eq:aw}
|\mathcal{S}f(z)| + e^{2\sigma(z)} |K(\widetilde{f}(z))| \leq \frac{2t}{(1-|z|^2)^2}\,,
\quad z\in\Bbb D
\end{equation}
for some $t<1$ and that
\begin{equation} \label{eq:sigma-condition-2}
\|\nabla\sg(z)\| \le \frac{C}{1-|z|^2},\quad z \in \D,
\end{equation}
for some $C <\infty$. Under these conditions we will show that the reflection $w^*=R(w)$ is quasiconformal.

Necessarily the analysis shifts to $\Sigma$ and some of the geometric notions attached to $\Sigma$ as a surface in $\mathbb{R}^3$ with its induced Euclidean metric $\mathbf{g}_0$, e.g., the gradient and the Hessian of a function, the covariant derivative and second fundamental form, and the curvature. As a reference we cite \cite{oneill:semiriemannian}, whose notation we generally follow. If $V$ is a vector field on $\Sigma$ we let $\overline{D}_V$ be the Euclidean covariant derivative on $\mathbb{R}^3$ in the direction $V$, applied to a function or a vector field on $\Sigma$, and we let $D_V$ be the covariant derivative on $\Sigma$. If $\psi$ is a function on $\Sigma$ then $\overline{D}_V \psi = D_V \psi = V\psi$. The gradient of $\psi$ is the vector field defined by
\[
\langle \nabla \psi, V \rangle = V\psi
\]
and its Hessian is the symmetric, covariant 2-tensor defined by
\[
\Hess\psi(V,W) = \langle D_V\nabla \psi ,W\rangle.
\]
If $W$ is a vector field on $\Sigma$ then
\[
\overline{D}_VW = D_VW + \sff(V,W)
\]
where $\sff(V,W)$ is the second fundamental form of $\Sigma$.

We can regard $w \mapsto R(w)$ as a vector field on $\Sigma$ (
not tangent to $\Sigma$) and we will compute its covariant derivative $\overline{D}_VR$ in the direction of a vector $V$, $\|V\|=1$, tangent to $\Sigma$. At each $w\in \Sigma$ we seek upper and lower bounds
\[
m(w) \le \| \overline{D}_VR \| \le M(w),
\]
where $\sup_{w\in\Sigma} M(w)/m(w)$ is bounded by a quantity depending on $t$ and $C$.

To do this we must translate the inequality \eqref{eq:aw} to one for functions defined on the surface. This requires the full differential-geometric definition of the conformal Schwarzian as a symmetric, traceless 2-tensor, and uses  in particular a generalization of the chain rule \eqref{eq:chain-rule} for the Schwarzian. We refer to \cite{cdo:injective-lift} for the details as they are applicable here, and to \cite{os:Schwarzian} for a more general treatment.

Very briefly, the main points are these. For a function $\psi$ defined on a 2-dimensional Riemannian manifold $(M,\g)$ the Schwarzian tensor of $\psi$ is
\begin{equation} \label{eq:B}
B_\g(\psi) = \Hess_\g\psi - d\psi \otimes d\psi - \frac{1}{2}(\Delta_\g \psi - \|\nabla_\g \psi\|_\g^2)\g
\end{equation}
where the Hessian, Laplacian, gradient, and norm are taken with respect to a Riemannian metric $\g$. The final term is the trace of $ \Hess_\g\psi - d\psi \otimes d\psi$, so the full tensor is traceless. If $f$ is a conformal mapping with conformal factor $e^{2\psi}\g$ then, by definition,
\[
\Sch_\g f = B_\g(\psi).
\]
In the case of a harmonic map $f$ and its lift $\widetilde{f}\colon (\D, |dz|^2) \rightarrow (\Sigma, {\mathbf{g_0}})$, with conformal factor $\widetilde{f}^*({\mathbf{g}_0}) = e^{2\sg}|dz|^2$ as before, we have
\[
\Sch f = \Sch \widetilde{f} = B(\sigma),
\]
with respect to the Euclidean metric, i.e., computing the right-hand side produces $2(\sigma_{zz} - \sigma_z^2)$, which we took as the definition of the harmonic Schwarzian. Here, and below, when a quantity is calculated with respect to the Euclidean metric we drop the subscript $\mathbf{g}_0$.

The quantities defining $B_\g(\psi)$ which depend on the metric change in a not very complicated manner when the metric changes \emph{conformally}. This is the basis for a generalized chain rule. It reads, in one form,
\[
B_{\hat{\mathbf{g}}}(\psi - \rho) = B_\g(\psi) - B_\g(\rho), \quad \hat{\g} = e^{2\rho}\g,
\]
 and (equivalently) in terms of conformal mappings, say $(M_1,\g_1) \stackrel{h}\longrightarrow (M_2,\g_2) \stackrel{f}\longrightarrow (M_3, \g_3)$,
 \[
 \Sch_{\g_1}(f \circ h) = h^*(\Sch_{\g_2} f) + \Sch_{\g_1} h.
 \]
 From the last equation, if $f$ and $h$ are inverse to each other then $\Sch_{\g_1} h  = - h^*(\Sch_{\g_2} f)$.

Specializing to our case, but set up a little differently than before, we find the following.  Recall from \eqref{eq:lambda-Sigma} the  metric $\lambda_\Sigma^2\, {\mathbf{g}_0}$  with $\lambda_\Sigma \circ \widetilde{f} = e^{-\sg}\lambda_\D$. Consider $\widetilde{f} \colon (\D, \g) \rightarrow (\Sigma, {\mathbf{g}_0})$, $\g= \lambda_\D^2|dz|^2$, as a conformal mapping with conformal factor $e^{2\sigma}\lambda_\D^{-2}$. We take the Schwarzian tensor of $\widetilde{f}$ with respect to $\g$:
\[
\Sch_\g \widetilde{f} = B_\g(\sg - \log\lambda_\D).
\]
Similarly, if $\widetilde{h} = \widetilde{f}^{-1}$ then $\widetilde{h}\colon (\Sigma, \mathbf{g}_0) \rightarrow (\D, \g)$  is conformal with conformal factor $\lambda_\Sigma^2$. The Schwarzian tensor of $\widetilde{h}$ is with respect to the induced Euclidean metric on $\Sigma$ and
\[
\Sch \widetilde{h} = B(\log \lambda_\Sigma).
\]
From the formulas above,
\[
B(\log \lambda_\Sigma) = \Sch \widetilde{h} = -\widetilde{h}^*\Sch_\g \widetilde{f} = -\widetilde{h}^*(B_\g(\sg - \log\lambda_\D)).
\]
while
\[
\begin{aligned}
B_\g(\sg - \log \lambda_\D) & = B(\sg) - B( \log \lambda_\D) = B(\sg),
\end{aligned}
\]
the last equation holding because one has $B(\log\lambda_\D) = 0$ (computing in the Euclidean metric).

On the other hand, $h\colon (\Sigma , \widetilde{\g}) \rightarrow (\D, \g)$ is an isometry for $\widetilde{\g} = \lambda_\Sigma^2\,\mathbf{g}_0$, thus
\[
\|B(\log \lambda_\Sigma) \|_{\widetilde{\g}}= \|B_\g(\sigma - \log\lambda_\D)\|_\g,
\]
and in turn
\[
 \|B_\g(\sigma - \log\lambda_\D)\|_\g =\|B(\sigma)\|_\g = \lambda_\D^{-2}\|B(\sigma)\|= \lambda_\D^{-2}|\Sch f|.
\]
In the final term $\Sch f$ is the harmonic Schwarzian. Combining these with \eqref{eq:aw} we find
\[
\begin{aligned}
\|B(\log \lambda_\Sigma) \|_{\widetilde{\g}} + \lambda_\Sigma^{-2}|K|&= \|B(\log \lambda_\Sigma) \|_{\widetilde{\g}} + \lambda_\D^{-2}e^{2\sg}|K|\\
& = \lambda_\D^{-2}(|\Sch f| +e^{2\sg}|K|) \le 2t.
\end{aligned}
\]
Finally, we switch to the norm in the Euclidean metric and state the results of the calculations above as a lemma.
\begin{lemma} \label{lemma:aw-on-Sigma}
If $f$ satisfies \eqref{eq:aw} then
\begin{equation} \label{eq:aw-Sigma}
\|B(\log \lambda_\Sigma)\| +|K| \le 2t\lambda_\Sigma^{2}.
\end{equation}
\end{lemma}
\noindent This is the inequality we use when working on $\Sigma$, eliminating direct mention of $\widetilde{f}$.

\smallskip

We proceed with the computation of $\overline{D}_VR$ using the formula \eqref{eq:R-gradient},
\[
R = {\rm Id}  +2J(\nabla\log\lambda_\Sigma),
\]
and the formula \eqref{eq:DJ},
\[
J'(x) = \frac{1}{\|x\|^4}(\|x\|^2{\rm Id}-2Q(x)).
\]
We have, first,
\[
\overline{D}_VR= V + 2J'(\nabla\log\lambda_\Sigma)(\overline{D}_V\nabla\log\lambda_\Sigma),
\]
and also the relation
\[
\overline{D}_V\nabla\log\lambda_\Sigma = D_V\nabla\log\lambda_\Sigma + \sff(V,\nabla\log\lambda_\Sigma).
\]
Hence
\[
\begin{aligned}
\overline{D}_VR & = V + \frac{2}{\|\nabla\log\lambda_\Sigma\|^4}\left\{\|\nabla\log\lambda_\Sigma\|^2{\rm Id} -2Q(\nabla\log\lambda_\Sigma)(D_V\nabla\log\lambda_\Sigma +\sff(V,\nabla\log\lambda_\Sigma))\right\}\\
\end{aligned}
\]
At this point it is prudent to simplify the notation somewhat. Let
\[
\Lambda = \|\nabla\log \lambda_\Sigma\|, \quad Q = Q(\nabla\log\lambda_\Sigma), \quad \sff = \sff(V,\nabla\log\lambda_\Sigma).
\]
Furthermore,
\[
\Hess(\log\lambda_\Sigma)(V,W) = \langle D_V\nabla\log\lambda_\Sigma, W\rangle
\]
so we identify the vector $D_V\nabla\log\lambda_\Sigma$ with the 1-tensor  $\Hess(\log\lambda_\Sigma)(V, \,\cdot\,)$
and write
\[
H= D_V\nabla\log\lambda_\Sigma.
\]
The Schwarzian tensor enters 
through the Hessian terms, but this is not immediate.

The expression for $\overline{D}_VR$ now appears a little more manageable:
\[
\overline{D}_VR=V+\frac{2}{\Lambda^2}\left\{H-\frac{2}{\Lambda^2}Q(H)+\sff-\frac{2}{\Lambda^2}Q(\sff)\right\}.
\]
To be clear, $\Lambda$ is a scalar, $\sff$ and $H$ are vectors,  and $Q$ is a matrix operating on the vectors $\sff$ and $H$.

To find the norm $\|\overline{D}_VR\|^2$  we are aided by several facts. First, $H$ is tangent to $\Sigma$ while $\sff$ is normal to $\Sigma$. Second, $Q$ is symmetric and
\[
Q^2 = \Lambda^2Q.
\]
 Finally, from its definition,
 \[
 Q_{ij} =  Q(\nabla\log\lambda_\Sigma)_{ij}= (\nabla\log\lambda_\Sigma)_i(\nabla\log\lambda_\Sigma)_j
 \]
 and it is easy to see that for any vector $X$ one has
\[
 Q(X) = \langle \nabla\log\lambda_\Sigma,X\rangle \nabla\log\lambda_\Sigma.
 \]
  Hence
\[
\langle Q(\sff),V\rangle = \langle \sff, Q(V)\rangle = \langle \nabla\log\lambda_\Sigma,V\rangle \langle\sff, \nabla\log\lambda_\Sigma\rangle =0,
\]
because $\sff$ is normal to $\Sigma$ and so is orthogonal to $\nabla\log\lambda_\Sigma$. In expanding $\|\overline{D}_VR\|^2$ a number of terms then drop out and, at length,  we obtain
\begin{equation} \label{eq:DVR}
\|\overline{D}_VR\|^2 = 1+ \frac{4}{\Lambda^2}\langle H,V\rangle +\frac{4}{\Lambda^4}\left\{\|H\|^2-2\langle Q(H),V\rangle +\|\sff\|^2\right\}
\end{equation}
where we have also used $\|V\|=1$.

 Referring to the definition \eqref{eq:B} we have
\[
B(\log\lambda_\Sigma) = \Hess(\log\lambda_\Sigma)-d\log\lambda_\Sigma \otimes d\log\lambda_\Sigma -\frac{1}{2}(\Delta \log\lambda_\Sigma - \|\nabla \log\lambda_\Sigma\|^2){\mathbf{g}_0}.
\]
Evaluate $B(\log\lambda_\Sigma)(V, \,\cdot \,)$ and treat this 1-tensor as a vector, which, continuing the pattern of notation, we will denote by $B$. With these abbreviations note that \eqref{eq:aw-Sigma} implies
\begin{equation} \label{eq:aw-Sigma-abbreviated}
\|B\|+ |K| \le 2t\lambda_\Sigma^2.
\end{equation}
Next, in components the 2-tensor $d\log\lambda_\Sigma \otimes d\log\lambda_\Sigma$ is exactly $Q(\nabla\log\lambda_\Sigma)$, which we have denoted by $Q$.  Finally we write
\[
\rho = \frac{1}{2}(\Delta \log\lambda_\Sigma - \|\nabla \log\lambda_\Sigma\|^2)= \frac{1}{2}(\Delta\log\lambda_\Sigma - \Lambda^2).
\]
for the trace. In these terms
\[
H = B + Q(V) + \rho V.
\]
and in \eqref{eq:DVR},
\[
\langle H, V\rangle = \langle B,V\rangle + \langle Q(V),V\rangle + \rho,
\]
\vskip1mm
\[
\|H\|^2 = \|B\|^2 + \Lambda^2\langle Q(V),V\rangle + \rho^2 + 2 \langle B, Q(V)\rangle +2\rho\langle B, V\rangle +2\rho\langle Q(V),V\rangle,
\]
\[
\langle Q(H),V\rangle=\langle H,Q(V)\rangle = \langle B,Q(V)\rangle +\Lambda^2\langle Q(V),V\rangle +\rho\langle Q(V),V\rangle.
\]

Substitution results in a quite compact expression:
\[
\|\overline{D}_VR\|^2= \frac{4}{\Lambda^4}\left\{\|B + \frac{1}{2}(\Delta \log\lambda_\Sigma)V\|^2+\|\sff\|^2\right\}.
\]
This is the penultimate form. The final step, to bring in the inequality \eqref{eq:aw-Sigma-abbreviated} for the Schwarzian,  is to introduce the curvature.

The curvature of $\Sigma$ with the metric $\lambda_\Sigma^2{\mathbf{g}_0}$ is $-4$ since  $(\Sigma, \lambda_\Sigma^2\,\mathbf{g}_0)$ is isometric to $(\D, \lambda_\D|dz|^2)$.
For the curvature $K \le 0$ of $\Sigma$ as a minimal surface one obtains
\[
\Delta \log\lambda_\Sigma = 4\lambda_\Sigma^2 - |K|.
\]
Hence
\begin{equation} \label{eq:DVR-final}
\|D_VR\|^2 = \frac{4}{\Lambda^4}\left\{\|B-\frac{1}{2}|K|V+2\lambda_\Sigma^2V\|^2 + \|\sff\|^2\right\}.
\end{equation}
We want to bound this from above and below.

\smallskip

To obtain a lower bound we drop the term $\|\sff\|^2$ and use \eqref{eq:aw-Sigma-abbreviated}:
\[
\begin{aligned}
\|D_VR\| &\ge \frac{2}{\Lambda^2}\|B-\frac{1}{2}|K|V+2\lambda_\Sigma^2V\|\ge \frac{2}{\Lambda^2}\left\{2\lambda_\Sigma^2 - \|-B+\frac{1}{2}|K|V\|\right\}\\
&\ge \frac{2}{\Lambda^2}\left\{2\lambda_\Sigma^2-\|B\|-\frac{1}{2}|K|\right\} \ge \frac{4\lambda_\Sigma^2}{\Lambda^2}(1-t).
\end{aligned}
\]

To obtain  an upper bound we have to estimate the term $\|\sff\|$. On a minimal surface we always have $\sff(X,Y) \le \sqrt{|K|}\|X\|\,\|Y\|$, and so for our case
\[
\|\sff\|=\|\sff(V,\nabla\log\lambda_\Sigma)\| \le \sqrt{|K|} \,\|\nabla\log\lambda_\Sigma\| = \sqrt{|K|}\Lambda.
\]
We need estimates for each of the factors on the right,  and this is where we use the assumption \eqref{eq:sigma-condition-2}, that
\[
\|\nabla\sg(z)\| \le \frac{C}{1-|z|^2}.
\]

An inequality for the curvature follows simply from dropping the positive $\|B\|$ term in \eqref{eq:aw-Sigma-abbreviated}, giving
\[
|K|\le 2t\lambda_\Sigma^2.
\]
Next, from $\log(\lambda_\Sigma\circ \widetilde{f}) = \log\lambda_\D - \sg$ and the bound on $\|\nabla\sg\|$ we have
\[
\begin{aligned}
e^{\sg(z)}\Lambda &= e^{\sg(z)}\|\nabla\log\lambda_\Sigma(\widetilde{f}(z))\| =\|\nabla\lambda_\D(z) - \nabla\sg(z)\|\\
& \le \|\nabla\log\lambda_\D(z)\| + \|\nabla\sigma(z)\| \le \frac{2+C}{1-|z|^2}.
\end{aligned}
\]
Multiplying through by $e^{-\sg}$ brings back $\lambda_\Sigma$ on the right:
\[
\Lambda \le (2+C)\lambda_\Sigma.
\]
Finally,
\[
\|\sff\|^2 \le |K|\Lambda^2 \le |K|(2+C)^2\lambda_\Sigma^2 \le2t(2+C)^2\lambda_\Sigma^4.
\]

Back to the equation \eqref{eq:DVR-final} for $\|D_VR\|^2$, we have
\[
\begin{aligned}
\|D_VR\| &\le \frac{2}{\Lambda^2}\left\{\|B-\frac{1}{2}|K|V+2\lambda_\Sigma^2V\|+\|\sff\|\right\}\\
&\le \frac{2}{\Lambda^2}\left\{\|B\| + \frac{1}{2}|K| + 2\lambda_\Sigma^2 +\|\sff\|\right\}\\
&\le \frac{2}{\Lambda^2}\left\{2t\lambda_\Sigma^2+2\lambda_\Sigma^2 +\sqrt{2t}(2+C)\lambda_\Sigma^2\right\}\\
&=\frac{2\lambda_\Sigma^2}{\Lambda^2}\left\{2t+\sqrt{2t}(2+C)+2\right\}.
\end{aligned}
\]
Combining the upper and lower bounds for $\|D_VR\|$ gives
\begin{equation} \label{eq:R-is-qc}
\frac{\max_{\|V\|=1}\|D_VR\|}{\min_{\|V\|=1}\|D_VR\|}\le \frac{2t+\sqrt{2t}(1+C) +2}{2(1-t)}.
\end{equation}
This shows that $R$ is quasiconformal as a mapping from $\Sigma$ to its reflection $\Sigma^*$. The extension of $\widetilde{f}$ to a mapping $\widetilde{F}: \overline{\mathbb{C}} \longrightarrow \overline{\Sigma} \cup \Sigma^*$ is as in \eqref{eq:F}. It, too, is quasiconformal with the same bound for the distortion. This completes the proof of Theorem \ref{theorem:extension}.

\bigskip

When $f$ is analytic satisfying the Ahlfors-Weill condition the quasiconformality of the reflection is measured simply by the Beltrami coefficient, and this turns out to be
\[
\mu(1/\bar{z}) = \frac{\partial_{\bar{z}} R(z)}{\partial_z R(z)} = -\frac{1}{2}(1-|z|^2)^2Sf(z), \quad z \in \D.
\]
Thus $|\mu| \le t <1$ and the extension of $f$ is a $(1+t)/(1-t)$-quasiconformal mapping of $\overline{\mathbb{C}}$. In the general case it is a question what one might take as a substitute for the Beltrami coefficient, but specializing to the analytic, planar case the bound \eqref{eq:R-is-qc} becomes
\[
\frac{\max_{\|V\|=1}\|D_VR\|}{\min_{\|V\|=1}\|D_VR\|}\le \frac{2t +2}{2(1-t)} = \frac{1+t}{1-t}
\]
because all estimates involving the curvature and the second fundamental form (and the upper bound for $\|\nabla\sg\|$) need not enter at all.

\subsection*{Remark:}  We have one final comment on when the condition
\[
\|\nabla \sg(z)\| \le \frac{C}{1-|z|^2}
\]
is satisfied if $f$ satisfies the injectivity condition \eqref{eq:nehari-harmonic}.

Suppose $u_{\widetilde{f}}$ has a critical point at $0$. This means that $\sigma_z(0)=0$ and we claim that
\begin{equation} \label{eq:new-sg-condition}
\|\nabla\sigma(z) \| \le \frac{2|z|}{1-|z|^2}.
\end{equation}
By applying a rotation of the disk it suffices to establish this on $[0,1)$. With $k(x) = \sigma_z(x)$ we find
\[
k'(x) = \sg_{zz}(x) + \sg_{z\bar{z}}(x)= (\sg_{zz}(x)-\sg_z(x)^2)+\sg_{z\bar{z}}(x)+k(x)^2,
\]
The bound  \eqref{eq:nehari-harmonic} says that
\[
 |\sg_{zz}(x)-\sg_z(x)^2|+2|\sg_{z\bar{z}}(x)| \le \frac{1}{(1-|z|^2)^2},
 \]
whence
\[
|k'(x)|  \le \frac{1}{(1-x^2)^2} + |k(x)^2|.
\]
Now let $a(x) = |k(x)|$, $b(x) = x/(1-x^2)$. Then
\[
a'(x) \le |h'(x)| \le \frac{1}{(1-x^2)^2} + a(x)^2 \quad \text{while} \quad b'(x) = \frac{1}{(1-x^2)^2} + b(x)^2.
\]
A standard comparison argument gives $a(x) \le b(x)$, which is our claim.

Suppose that the surface $\Sigma$ is bounded, or equivalently that $u_{\widetilde{f}}$ has a critical point somewhere in the disk. We may compose with a M\"obius transformation of $\D$ onto itself to locate the critical point at the origin, and for the new conformal factor we will have \eqref{eq:new-sg-condition}. Since the new and original conformal factors are scaled by a factor that is smooth in the closed disk, \eqref{eq:new-sg-condition} will also hold for the original conformal factor up to a constant multiple.

\bibliographystyle{amsplain}

\bibliography{extension}

\end{document}